\documentclass[12pt]{article}

\usepackage{amsthm,color,amsmath,amsfonts,amssymb}

\parskip=\medskipamount
\def\onehalf{{\textstyle\frac12}}

\def\jet{J^1\tau}
\def\dualJ{J^1\tau^*}

\def\cinfty#1{C^{\scriptscriptstyle\infty}(#1)}
\def\vectorfields#1{{\mathcal X}(#1)}
\def\Vvectors{{\mathcal X}_V(E)}
\def\oneforms#1{{\mathcal X}^*(#1)}
\def\forms#1#2{{\textstyle\bigwedge}^{#1}(#2)}

\def\tvectors{{\mathcal X}_t(E)}

\def\lie#1{{\mathcal L}_{#1}}
\def\fpd#1#2{\frac{\partial #1}{\partial #2}}
\def\R{\mathbb{R}}
\def\D#1{{\mathcal D}_{(#1)}}

\def\proof{{\sc Proof.}}

\newtheorem{prop}{\bf Proposition}
\newtheorem{lemma}{\bf Lemma}
\newtheorem{thm}{\bf Theorem}
\newtheorem{dfn}{\bf Definition}

\def\h#1{{#1}^{h}}
\def\V#1{{#1}^{v}}
\def\wt#1{\widetilde{#1}}

\def\sp{\mathop{\rm sp}}

\def\hook2{\mathop{\kern -2pt\hbox to 6pt{\hrulefill}
                      \hbox{\vrule\phantom{\vbox to 8pt{}}}}\kern 1pt}

\begin{document}

\title{Lifting geometric objects to the dual of the first jet bundle of a bundle fibred over $\R$}

\author{W.\ Sarlet$^{a,b}$ and G.\ Waeyaert$^{a}$\\
{\small ${}^a$Department of Mathematics, Ghent University }\\
{\small Krijgslaan 281, B-9000 Ghent, Belgium}\\[1mm]
{\small ${}^b$Department of Mathematics and Statistics, La Trobe University}\\
{\small Bundoora, Victoria 3086, Australia}
}
\date{}
\maketitle

\begin{quote}
{\bf Abstract.} {\small We study natural lifting operations from a bundle $\tau:E\rightarrow\R$ to the bundle $\pi:\dualJ\rightarrow E$ which is the dual of the first-jet bundle $\jet$. The main purpose is to define a complete lift of a type $(1,1)$ tensor field on $E$ and to understand all features of its construction. Various other lifting operations of tensorial objects on $E$ are needed for that purpose. We prove that the complete lift of a type $(1,1)$ tensor with vanishing Nijenhuis torsion gives rise to a Poisson-Nijenhuis structure on $\dualJ$, and discuss in detail how the construction of associated Darboux-Nijenhuis coordinates can be carried out.}
\end{quote}

\section{Introduction}
Consider a bundle $\tau:E\rightarrow \R$ with $\dim E=n+1$. Local coordinates on $E$ will be denoted by $(t,q^i)$. The construction of the first jet manifold $\jet$ of $\tau$ is well known (see e.g.\ \cite{Saunders}). It is an affine bundle over $E$ with induced coordinates $(t,q,\dot{q})$ say, and can be considered as an affine sub-bundle of $TE$, locally determined by coordinates of the form $(t,q,\dot{t}=1,\dot{q})$. $\jet$ is the space to be for the analysis of intrinsic or geometric aspects of time-dependent Lagrangian systems. The motivation for the present paper is more like setting the stage for future applications on the Hamiltonian side. So, we consider the cotangent bundle $T^*E$, with natural coordinates $(t,q^i,p_0,p_i)$, which is said to be the extended dual $(\jet)^\dagger$ of $\jet$. The main space of interest, however, is the \emph{dual of $\jet$}, denoted by $\dualJ$, which is the quotient space $T^*E/\langle dt\rangle$ and is sometimes called the vertical c
 otangent
  bundle. There are natural projections $\rho:T^*E\rightarrow \dualJ$ and $\pi:\dualJ\rightarrow E$. Each point $m\in\dualJ$ is an equivalence class of covectors $\langle\alpha\rangle\!\!\!\mod dt$ at $\pi(m)$; saying that $m$ has coordinates $(t,q^i,p_i)$ means that a representative of the class is given by $\alpha_{(t,q)}=p_idq^i$. Elements of $\dualJ$ have a well-defined action on vertical tangent vectors to $E$:
\[
\mbox{for}\ \ v= \left.v^i\fpd{}{q^i}\right|_{(t,q)} \in T_{(t,q)}E,\ \mbox{ we have\ } \langle v, \langle\alpha\rangle \rangle = v^ip_i.
\]
$\dualJ$ is the natural bundle for the description of time-dependent Hamiltonian systems: a Hamiltonian is a section $h$ of the line bundle $\rho:T^*E\rightarrow\dualJ$. Locally, $h$ defines a function $H$ on $\dualJ$, determined by $h:(t,q,p)\mapsto (t,q,p_0=-H(t,q,p),p)$ say (where the minus sign is of course a matter of convention). A Hamiltonian system then is a vector field $X_h$ on $\dualJ$, satisfying $i_{X_h}h^*\omega_E=0, \langle X_h,dt\rangle=1$, where $\omega_E$ is the canonical symplectic form on $T^*E$, so that locally $h^*\omega_E = dp_i\wedge dq^i - dH\wedge dt$.

An important point to note here is that in many papers on time-dependent mechanical systems, $E$ is identified with a trivial bundle $\R\times M$, which subsequently implies that $\jet\cong \R\times TM$ and $\dualJ\cong \R\times T^*M$. The list of citations we could insert here is endless, but to name just a few, in different contexts, see \cite{ACI}, \cite{CaRan}, \cite{CMS}, \cite{EMR}. Often, the trivial bundle $\R\times M$ is chosen as model for the configuration space from the outset. A more sophisticated argument, however, is that a bundle such as $\tau: E\rightarrow \R$ can be identified with $\R\times M$ by choosing a trivialization. So this is fine, as long as one keeps in mind that the essence of applications in time-dependent mechanics is to allow for time-dependent coordinate transformations, which are transformations which do not preserve the product bundle structure $\R\times M$. Our goal specifically in this paper is to study natural lifting operations from ten
 sor fiel
 ds on $E$ to corresponding tensor fields on $\dualJ$. Therefore, we want to avoid making such identifications because this would hold the danger of introducing lifting operations which are natural only under transformations which preserve the corresponding product structure (examples of such operations can be found e.g.\ in \cite{Vondra}). In this respect, our approach more closely relates to the work of Sardanashvily and co-workers on time-dependent mechanics, see \cite{Sarda} and \cite{GMS}. Other similar settings can be found e.g.\ in \cite{MV}.

Starting from objects on $E$, corresponding lifted objects on $TE$, $T^*E$ or $\jet$ are rather well known. This is less the case for lifting operations to $\dualJ$, with the exception of \cite{GMS} already cited. Our main objective is to come to an intrinsic definition of the \emph{complete lift of a type $(1,1)$ tensor field from $E$ to $\dualJ$} and to study its properties. To the best of our knowledge, this has not been analysed before and earlier work of one of us \cite{CCS} about lifting to a cotangent bundle will serve as a source of inspiration. For some other interesting background information about general aspects of natural operations, see \cite{KMS}. Another standard reference for lifting operations is \cite{YI}.

In the next section, we discuss various lifts of vector fields and 1-forms to $\dualJ$ and list some immediate properties. Ways of lifting type $(1,1)$ tensors from $E$ to $\dualJ$ are introduced in Section~3. Further properties relating the constructions of the two preceding sections are derived in Section~4. They are indispensable for proving in Section~5 that the canonical Poisson structure on $\dualJ$, together with the complete lift of a type $(1,1)$ tensor $R$ on $E$ with vanishing Nijenhuis torsion, determine a Poisson-Nijenhuis structure on $\dualJ$. The construction of Darboux-Nijenhuis coordinates for this structure is explained in detail in Section~6. In the final section, we say a few words about the applications to time-dependent Hamiltonian systems we have in mind for future publications.

\section{Lifting vector fields and 1-forms}

Consider $\Vvectors$, the $\cinfty{E}$-module of vertical vector fields on $E$. For each $X\in\Vvectors$ we can define a fibre linear function on $\dualJ$ as follows.

\begin{dfn} If $X\in\Vvectors$, we denote by $F_X\in\cinfty{\dualJ}$ the function which at each point $m\in\dualJ$ takes the value $F_X(m) = \langle X_{\pi(m)},m\rangle$.
\end{dfn}
If $m$ has coordinates $(t,q,p)$ and $X=X^i(t,q)\partial/\partial q^i$, the coordinate expression of $F_X$ is given by $F_X(t,q,p) = p_iX^i(t,q)$.

A first set of interesting vector fields on $\dualJ$ is obtained by vertically lifting 1-forms on $E$.

\begin{dfn} For $\alpha\in\oneforms{E}$, $\V{\alpha}\in\vectorfields{\dualJ}$ is determined by
\begin{align*}
\V{\alpha}(\pi^*f) &=0, \quad \forall f\in\cinfty{E} \\
\V{\alpha}(F_X) &= \pi^*\langle X,\alpha\rangle \quad \forall X\in\Vvectors.
\end{align*}
$\V{\alpha}$ is called the vertical lift of $\alpha$.
\end{dfn}
In coordinates, if $\alpha=\alpha_0(t,q) dt + \alpha_i(t,q) dq^i$, we have
\begin{equation}
\V{\alpha}=\alpha_i\fpd{}{p_i}.  \label{valpha}
\end{equation}
It is worth observing that in fact the assignment $\alpha \mapsto \V{\alpha}$ is a map from $\oneforms{E}/\langle dt\rangle$ into $\vectorfields{\dualJ}$.

$\dualJ$ does not carry a canonical 1-form, but there is a canonical equivalence class of 1-forms modulo the module generated by $dt$, which we shall denote by $\langle\theta\rangle$. As an element of $\oneforms{\dualJ}/\langle dt\rangle$, it is meant to have a well-defined action, at each point $m\in\dualJ$, on vectors which are vertical with respect to the projection $\tau\circ\pi:\dualJ \rightarrow \R$.

\begin{dfn} The equivalence class $\langle\theta\rangle\in\oneforms{\dualJ}/\langle dt\rangle$ is defined as follows: $\forall m\in\dualJ$ and $v_m\in T_m(\dualJ)$, vertical with respect to $\tau\circ\pi$, we have
\[
\langle v_m, \langle\theta\rangle_m\rangle = \langle T\pi(v_m),m\rangle.
\]
\end{dfn}
In coordinates, $\langle\theta\rangle= p_i dq^i\!\!\! \mod dt$. It follows that
\begin{equation}
\Theta = p_i dq^i \wedge dt \label{bigtheta}
\end{equation}
is a canonically defined 2-form on $\dualJ$. It can be characterized alternatively by the following property, which mimics a well known characterization of the canonical 1-form on a cotangent bundle. A section of the bundle $\pi:\dualJ\rightarrow E$ is an element $\langle\alpha\rangle\in\oneforms{E}/\langle dt\rangle$ and for each representative $\alpha$ of the class we have that $\alpha^*\Theta= \alpha\wedge dt$.

We now have the tools to define the complete lift of two classes of vector fields on $E$, one is the module of vertical vector fields $\Vvectors$ already introduced; the other is the set of vector fields with the property $\langle X,dt\rangle =1$ which we shall denote by $\tvectors$. Elements of $\tvectors$ can be regarded as sections of $\jet\rightarrow E$, if $\jet$ is seen as the submanifold of $TE$ described before.

\begin{dfn} For all $X\in\Vvectors \cup\tvectors$, the complete lift $\widetilde{X}\in\vectorfields{\dualJ}$ is defined by the following two requirements:
(i)\ $\widetilde{X}$ is $\pi$-related to $X$; (ii)\ $\lie{\widetilde{X}}\Theta =0$.
\end{dfn}
It is an easy computation to check that in coordinates,
\begin{align}
\widetilde{X} &= X^i \fpd{}{q^i} - p_j\fpd{X^j}{q^i}\fpd{}{p_i}, \quad \mbox{for}\ X=X^i(t,q)\fpd{}{q^i}\in\Vvectors \label{complete1} \\
\widetilde{X} &= \fpd{}{t} + X^i \fpd{}{q^i} - p_j\fpd{X^j}{q^i}\fpd{}{p_i}, \quad \mbox{for}\ X=\fpd{}{t} + X^i(t,q)\fpd{}{q^i}\in\tvectors. \label{complete2}
\end{align}
Both types of complete lifts are introduced in \cite{GMS} in a different way, mainly based on coordinate calculations. In any event, it is an instructive exercise to verify by direct coordinate calculations that all lifting constructions so far introduced are indeed well defined, i.e.\ behave properly under a time-dependent coordinate transformation $t=t,\ Q^i=Q^i(t,q)$ on $E$ and the induced transformation $(t,q,p)\mapsto (t,Q,P)$ on $\dualJ$, where
\[
P_j= p_i\fpd{q^i}{Q^j}(t,Q(t,q)).
\]
Some immediate properties of $\widetilde{X}$, which are easy to verify in coordinates, are
\begin{align*}
\mbox{for}\ X\in\Vvectors,& \quad i_{\widetilde{X}}\Theta = F_X dt, \\
\mbox{for}\ X\in\tvectors,& \quad i_{\widetilde{X}}\Theta\wedge dt = - \Theta.
\end{align*}

The main motivation for introducing the vector fields $\V{\alpha}$ and both types of $\widetilde{X}$, is that together they provide a local basis of vector fields on $\dualJ$. As such, they are perfectly suited to define other tensorial objects on $\dualJ$ in a coordinate free way, as we will see in the subsequent sections. It will then be interesting to have expressions for the Lie brackets of such vector fields. For all $\alpha,\beta\in\oneforms{E}$ and $X\in\Vvectors\cup\tvectors$, we have
\begin{align}
[\V{\alpha},\V{\beta}] &= 0, \label{vbracket} \\
[\widetilde X, \V{\alpha}] &= \V{(\lie{X}\alpha)}, \label{tildevbracket} \\
[\widetilde X, \widetilde Y] &= \widetilde{[X,Y]}. \label{tildebracket}
\end{align}

\section{Lifting type $(1,1)$ tensor fields}

In what follows, $R$ will always denote a type $(1,1)$ tensor field on $E$ with the property $R(dt)=0$. We make no notational distinction between the action of a $(1,1)$ tensor on vector fields and its adjoint action on 1-forms; for example, for $X\in\vectorfields{E}$ and $\alpha\in\oneforms{E}$, we have $\langle R(X),\alpha\rangle = \langle X, R(\alpha)\rangle$. In coordinates, the tensor fields under consideration have the form
\begin{equation}
R = R^i_j(t,q)\fpd{}{q^i}\otimes dq^j + R^i_0(t,q)\fpd{}{q^i}\otimes dt. \label{R}
\end{equation}

\begin{dfn} The vertical lift $\V{R}$ is a vector field on $\dualJ$, determined by
\begin{align*}
\V{R}(\pi^*f) &=0, \quad \forall f\in\cinfty{E} \\
\V{R}(F_X) &= F_{R(X)} \quad \forall X\in\Vvectors.
\end{align*}
\end{dfn}
In coordinates,
\begin{equation}
\V{R}= p_iR^i_j\fpd{}{p_j}. \label{vR}
\end{equation}

\begin{dfn} The horizontal lift $\h{R}$ is a 1-form on $\dualJ$, which pointwise is defined by $\h{R}_m = R_{\pi(m)}(m)$ for all $m\in\dualJ$.
\end{dfn}
Recall that $m$ is not a covector at $\pi(m)$, but an equivalence class of covectors $\!\!\!\mod dt$. But the action of $R_{\pi(m)}$ on such a class is well defined because of the property $R(dt)=0$. $\h{R}$ is a semi-basic 1-form on $\dualJ$, which in coordinates reads
\begin{equation}
\h{R} = p_iR^i_j dq^j + p_iR^i_0 dt. \label{hR}
\end{equation}

Having added a new type of vertical vector field to the picture, it is appropriate that we complement the Lie bracket properties listed at the end of the previous section. For $\alpha\in\oneforms{E}$ and $X\in\Vvectors\cup\tvectors$, we have
\begin{align}
[\V{\alpha},\V{R}] &= \V{R(\alpha)}, \label{vbracket2} \\
[\V{R}_1,\V{R}_2] &= \V{(R_1R_2-R_2R_1)}, \label{vbracket3} \\
[\widetilde{X},\V{R}] &= \V{(\lie{X}R)}. \label{tildevbracket2}
\end{align}
Note in passing that if $R(dt)=0$, then for all $X\in\Vvectors\cup\tvectors$, also $(\lie{X}R)(dt)=0$.

To define the complete lift of a $(1,1)$ tensor to the cotangent bundle, one normally makes use of the canonical symplectic form. We don't have a symplectic structure on $\dualJ$.  Therefore, we look at relations which were established as properties in the cotangent case (see \cite{CCS}) as a source of inspiration to come to an alternative definition here. The following lemma will be useful for that purpose.

\begin{lemma} If $f$ is an arbitrary function on $E$, we have
\begin{align*}
\V{(f\alpha)} &= f\,\V{\alpha}, \quad \alpha\in\oneforms{E} \\
\V{(fR)} &= f \,\V{R}, \quad R\ \mbox{$(1,1)$-tensor on $E$} \\
\wt{fX} &= f\,\wt{X} - F_X\,\V{(df)}, \quad X\in\Vvectors \\
\lie{fX}R &= f\,\lie{X}R - X\otimes R(df) + R(X)\otimes df.
\end{align*}
\end{lemma}
\proof\ Follows from a straightforward computation, for example in coordinates. For completeness: the factor $f$ on the right-hand side of the first three relations should actually be $\pi^*(f)$, but we try to avoid an overload of notations here and in what follows. \qed

\begin{thm} Given a type $(1,1)$ tensor field $R$ on $E$ with the property $R(dt)=0$, there is a unique type $(1,1)$ tensor $\wt{R}$ on $\dualJ$, which has the properties
\begin{align}
\wt{R}(\V{\alpha}) &= \V{R(\alpha)}, \quad \forall\,\alpha\in\oneforms{E} \label{tRv} \\
\wt{R}(\wt{X}) &= \wt{R(X)} + \V{(\lie{X}R)}, \quad\forall\, X\in\Vvectors\cup\tvectors. \label{tRtX}
\end{align}
$\wt{R}$ is called the complete lift of $R$ to $\dualJ$.
\end{thm}
\proof\ Note first that for both types of vector fields in (\ref{tRtX}), we have that $R(X)$ is vertical, so that $\wt{R(X)}$ makes sense. As indicated before, the set of vector fields considered in the above relations constitutes a local basis for the vector fields on $\dualJ$. To be specific, we need $n$ independent $\alpha\in\oneforms{E}$, $n$ independent $X\in\Vvectors$ and one $X\in\tvectors$ to generate such a basis. Imposing linearity over the module $\cinfty{\dualJ}$ then further fixes $\wt{R}$ for all vector fields in $\vectorfields{\dualJ}$. But for this approach to be consistent, we need to verify that our construction does not depend on a specific selection of 1-forms $\alpha$ and vector fields $X$ on $E$. Since every other selection of independent $\alpha$ and $X$ will originate from a linear combination over $\cinfty{E}$ of the ones first thought of, the issue is to check that the defining relations behave properly under $f$-linear changes of $\alpha\in\oneforms
 {E}$ and
  $X\in\Vvectors$, with $f\in\cinfty{E}$. Using the results of the preceding lemma, we have that
\[
\wt{R}(\V{(f\alpha)}) = \wt{R}(f\V{\alpha}) = f\,\wt{R}(\V{\alpha})= f\,\V{(R(\alpha))} = \V{(R(f\alpha))}.
\]
Also, on the one hand
\begin{align*}
\wt{R}(\wt{fX}) &= \wt{R}\big(f\wt{X} - F_X\V{(df)}\big)= f\wt{R}(\wt{X}) - F_X\wt{R}(\V{(df)}) \\
&= f\,\big(\wt{R(X)} + \V{(\lie{X}R)}\big) - F_X\V{(R(df))}.
\end{align*}
This can be seen to match the sum of the following two expressions:
\[
\wt{R(fX)} = \wt{fR(X)} = f\,\wt{R(X)} - F_{R(X)}\V{(df)},
\]
and
\[
\V{(\lie{fX}R)} = f\V{(\lie{X}R)} - F_X\V{(R(df))} + F_{R(X)}\V{(df)},
\]
where we have used also the general property $\V{(Y\otimes\beta)}= F_Y\,\V{\beta}$. \qed

In coordinates, the complete lift of the tensor $R$ in (\ref{R}) reads,
\begin{align}
\wt{R} &= R^i_j \left(\fpd{}{q^i}\otimes dq^j + \fpd{}{p_j}\otimes dp_i\right) + R^i_0 \fpd{}{q^i}\otimes dt \nonumber \\
& \mbox{} + p_i\left(\fpd{R^i_j}{q^k} - \fpd{R^i_k}{q^j}\right) \fpd{}{p_j}\otimes dq^k  + p_i\left(\fpd{R^i_k}{t} - \fpd{R^i_0}{q^k}\right) \fpd{}{p_k}\otimes dt. \label{tildeR}
\end{align}

At this point, it is of some interest to compare $\wt{R}$ with the complete lift of $R$ to the cotangent bundle $T^*E$, let us call it $\wt{R}_{T^*}$, which in coordinates $(t,q^i,p_0,p_i)$ on $T^*E$ is given by
\begin{align*}
\wt{R}_{T^*} &= R^i_j \left(\fpd{}{q^i}\otimes dq^j + \fpd{}{p_j}\otimes dp_i\right) + R^i_0 \left(\fpd{}{q^i}\otimes dt + \fpd{}{p_0}\otimes dp_i\right)\\
& \mbox{} + p_i\left(\fpd{R^i_j}{q^k} - \fpd{R^i_k}{q^j}\right) \fpd{}{p_j}\otimes dq^k  + p_i\left(\fpd{R^i_k}{t} - \fpd{R^i_0}{q^k}\right) \fpd{}{p_k}\otimes dt \\
& \mbox{} + p_i\left(\fpd{R^i_0}{q^k} - \fpd{R^i_k}{t}\right)\fpd{}{p_0}\otimes dq^k.
\end{align*}
Considering the projection $\rho:T^*E\rightarrow \dualJ$, the concept of $\rho$-related vector fields is well known.
\begin{dfn} Type $(1,1)$ tensor fields $U$ on $T^*E$ and $V$ on $\dualJ$ are said to be $\rho$-related, if for all $\rho$-related pairs of vector fields $(Y,Z)$, we have that $U(Y)$ is $\rho$-related to $V(Z)$.
\end{dfn}
The following alternative characterization is easy to prove.
\begin{prop} $U$ and $V$ are $\rho$-related  if for all $\sigma\in\oneforms{\dualJ}$ we have that $U(\rho^*\sigma) = \rho^*(V(\sigma))$.
\end{prop}
\begin{prop} $\wt{R}_{T^*}$ on $T^*E$ and $\wt{R}$ on $\dualJ$ are $\rho$-related.
\end{prop}
This is easy to see  from the coordinate expressions by appealing to the result of the preceding proposition.

It is important for the next sections that we also pin down the adjoint action of $\wt{R}$ by using a natural local basis of 1-forms on $\dualJ$. Such a natural basis is being provided by pull backs of 1-forms on $E$, complemented by 1-forms of the type $dF_X$, with $F_X$ as introduced in Definition~1. Again, a straightforward coordinate calculation is sufficient to verify that the adjoint action of $\wt{R}$ is determined by the following coordinate free relations, in which also the horizontal lift of $(1,1)$ tensors exhibits its relevance.
\begin{prop} The action on 1-forms of the complete lift $\wt{R}$ on $\dualJ$ is fully determined by the relations
\begin{align}
\wt{R}(\pi^*\alpha) &= \pi^* R(\alpha), \quad \forall\,\alpha\in\oneforms{E} \label{tRa} \\
\wt{R}(dF_X) &= dF_{R(X)} - \h{(\lie{X}R)}, \quad \forall\,X\in\Vvectors. \label{tRdF}
\end{align}
\end{prop}

\section{Further properties of $\wt{R}$}

The main goal for this section is to compute the Nijenhuis torsion of the complete lift $\wt{R}$. For that we will need some auxiliary properties, for example information about the Lie derivatives of $\wt{R}$ with respect to vector fields of type $\V{\alpha}$ or $\wt{X}$. In turn, this is prompting for one further lifting operation, from a 2-form on $E$ to a type $(1,1)$ tensor field on $\dualJ$.

\begin{dfn} For $\omega\in\forms{2}{E}$ we define a type $(1,1)$ tensor $\V{\omega}$ on $\dualJ$, called the vertical lift of $\omega$, by the relations
\begin{align}
\V{\omega}(\V{\alpha}) &= 0, \quad \forall\,\alpha\in\oneforms{E} \label{2v1} \\
\V{\omega}(\wt{X}) &= \V{(i_X\omega)}, \quad \forall\, X\in\Vvectors\cup\tvectors. \label{2v2}
\end{align}
\end{dfn}
The defining relations are obviously linear here with respect to multiplication of $\alpha$ or $X$ with a function on $E$, so that this is well-defined. In coordinates, if
\begin{equation}
\omega = \onehalf \omega_{ij}(t,q) dq^i\wedge dq^j + \omega_{0i}(t,q)dt\wedge dq^i, \label{omega}
\end{equation}
then
\begin{equation}
\V{\omega} = \omega_{ij}\fpd{}{p_j}\otimes dq^i + \omega_{0j}\fpd{}{p_j}\otimes dt. \label{vomega}
\end{equation}
\begin{prop} The basic Lie derivatives of the complete lift $\wt{R}$ have the following expressions
\begin{align}
\lie{\wt{X}}\wt{R} &= \wt{\lie{X}R}, \quad \forall\, X\in\Vvectors\cup\tvectors \label{lieR1} \\
\lie{\V{\alpha}}\wt{R} &= \V{\big(-R\hook2 d\alpha + d(R\alpha)\big)}, \quad \forall\,\alpha\in\oneforms{E}. \label{lieR2}
\end{align}
\end{prop}
\proof\ The proof is a simple matter of evaluating both sides of the above claims on a basis of vector fields on $\dualJ$. We leave the first one as an exercise for the reader and show how it works for (\ref{lieR2}). We have
\[
\lie{\V{\alpha}}\wt{R}(\V{\beta}) = \lie{\V{\alpha}}\big(\wt{R}(\V{\beta})\big) - \wt{R}([\V{\alpha},\V{\beta}]) = 0,
\]
while
\begin{align}
\lie{\V{\alpha}}\wt{R}(\wt{X}) &= \lie{\V{\alpha}}\big(\wt{R}(\wt{X})\big) - \wt{R}([\V{\alpha},\wt{X}]) \nonumber \\
&= \lie{\V{\alpha}}\big(\wt{RX} + \V{(\lie{X}R)}\big) + \wt{R}\big(\V{(\lie{X}\alpha)}\big) \nonumber \\
&= - \V{(\lie{RX}\alpha)} + \V{\big(\lie{X}R(\alpha)\big)} + \V{\big(R(\lie{X}\alpha)\big)} \nonumber \\
&= \V{\big(-\lie{RX}\alpha + \lie{X}(R\alpha)\big)}  \label{aux} \\
&= \V{\big( -i_{RX}d\alpha + i_Xd(R\alpha)\big)} \nonumber \\
&= \V{\big(-i_X(R\hook2 d\alpha) + i_Xd(R\alpha)\big)}, \nonumber
\end{align}
from which the equality (\ref{lieR2}) now readily follows. In making this computation, we have made use of properties such as (\ref{vbracket}), (\ref{tildevbracket}) and (\ref{vbracket2}) and of course the defining relations of Theorem~1. \qed

Note in passing that for a 2-form $\omega$, $R\hook2\omega$ is not the same as $i_R\omega$; by $R\hook2\omega$ we mean the 2-form defined by $R\hook2\omega(X,Y)=\omega(RX,Y)$.

Now recall that the Nijenhuis torsion of a tensor such as $R$ is the type $(1,2)$ tensor $N_R$ defined by
\begin{equation}
N_R(X,Y) = [RX,RY] + R^2([X,Y]) - R([RX,Y]) - R([X,RY]),  \label{NR}
\end{equation}
with $X,Y\in\vectorfields{E}$. We also put $(i_X N_R)(Y) = N_R(X,Y)$ and recall then the property that for the action on vector fields
\begin{equation}
i_X N_R = \lie{RX}R - R\circ\lie{X}R. \label{XhookNR}
\end{equation}
For the computation of the Nijenhuis torsion of $\wt{R}$ we will need the following result.
\begin{prop} The action of $\wt{R}^2$ on vector fields on $\dualJ$ is determined by
\begin{align*}
\wt{R}^2(\V{\alpha}) &= \V{(R^2(\alpha))}, \quad \alpha\in\oneforms{E} \\
\wt{R}^2(\wt{X}) &= \wt{R^2(X)} + \V{(\lie{X}R^2)} + \V{(i_X N_R)}, \quad X\in\Vvectors\cup\tvectors.
\end{align*}
\end{prop}
\proof\ Follows easily from the defining relations in Theorem~1. \qed

\begin{prop} The Nijenhuis torsion of the complete lift $\wt{R}$ on $\dualJ$ is determined by the following relations, where $\alpha,\beta\in\oneforms{E}$ and $X,Y\in\Vvectors\cup\tvectors$,
\begin{align*}
N_{\wt{R}}(\V{\alpha},\V{\beta}) &= 0, \\
N_{\wt{R}}(\wt{X},\V{\alpha}) &= \V{((i_X N_R)(\alpha))}, \\
N_{\wt{R}}(\wt{X},\wt{Y}) &= \wt{N_R(X,Y)} + \V{(i_{[X,Y]} N_R)} + \V{\big(\lie{Y}(i_X N_R) - \lie{X}(i_Y N_R)\big)}.
\end{align*}
\end{prop}
\proof\ The proof is a matter of making use of the defining relations in Theorem~1 again, together with the results about $\wt{R}^2$ of the preceding proposition and a number of the bracket relations established before. But formally, all these relations are identical to the ones we know from the cotangent bundle situation. Specifically, what are defining relations here, namely (\ref{tRv},\ref{tRtX}), were properties of the complete lift to a cotangent bundle in \cite{CCS} (cf.\ Theorem~1 in that paper). Hence, we can simply refer to the proof of the Nijenhuis properties in \cite{CCS}. Since that proof was given in great detail for $N_{\wt{R}}(\wt{X},\wt{Y})$, we limit ourselves here to giving a sketch of the calculation for $N_{\wt{R}}(\wt{X},\V{\alpha})$. Starting from the definition of $N_{\wt{R}}$ and a first implementation of known properties, mainly from Theorem~1, we have
\begin{align*}
N_{\wt{R}}(\wt{X},\V{\alpha}) &= [\wt{RX},\V{(R\alpha)}] + [\V{(\lie{X}R)},\V{(R\alpha)}] + \wt{R}^2(\V{(\lie{X}\alpha)}) \\
& \ \ - \wt{R}([\wt{RX},\V{\alpha}]) - \wt{R}([\V{(\lie{X}R)},\V{\alpha}]) - \wt{R}([\wt{X},\V{(R\alpha)}]).
\end{align*}
A subsequent use of known bracket relations reduces the right-hand side to the vertical lift of the following aggregation of terms:
\[
\lie{RX}(R\alpha) - \lie{X}R(R\alpha) + R^2(\lie{X}\alpha) - R(\lie{RX}\alpha) + R(\lie{X}R(\alpha)) - R(\lie{X}(R\alpha)).
\]
It is now a simple matter to simplify this expression further to
\[
\lie{RX}R(\alpha) - \lie{X}R(R\alpha) = (i_X N_R)(\alpha),
\]
where one has to keep in mind that the order of composition of $(1,1)$ tensors changes, when passing from the action on vector fields to the adjoint action on 1-forms. \qed

\begin{thm} $N_{\wt{R}}=0$ if and only if $N_R=0$.
\end{thm}
\proof\ From the results of the preceding proposition, it is clear that $N_R=0$ implies $N_{\wt{R}}=0$. Conversely, $N_{\wt{R}}=0$ implies in particular that $\V{((i_X N_R)(\alpha))}=0$ for all 1-forms $\alpha$ on $E$ and $X\in\Vvectors\cup\tvectors$, which is equivalent to saying that the 1-form $(i_X N_R)(\alpha)=0 \!\mod dt$. In turn, this can be expressed as $\langle Y,(i_X N_R)(\alpha)\rangle=0$ for all vertical $Y$. Looking at $N_R$ again as vector-valued two-form, the conclusion is that
\[
N_R(X,Y)=0, \quad \forall\,Y\in\Vvectors,\  X\in\Vvectors\cup\tvectors.
\]
Thinking for a moment in coordinates, if we take an $X$ in $\tvectors$, we will still get zero by adding to the vertical $Y$ a term spanned by $\partial/\partial t$ because of the skew-symmetry of $N_R$. In the end, the fact that $N_R$ is known to be tensorial guarantees that $N_R(X,Y)=0$ for all $X,Y\in\vectorfields{E}$. \qed

\section{A Poisson-Nijenhuis structure on $\dualJ$}

It is known (see e.g.\ \cite{GMS}) that $\dualJ$ carries a canonical Poisson structure, namely the structure inherited from the Poisson structure on $T^*E$ via the projection $\rho:T^*E\rightarrow \dualJ$. The point is that the Poisson bracket of functions $\rho^*F$ and $\rho^*G$, with $F,G\in\cinfty{\dualJ}$ is a function of the same type, constant on fibres. According to a general result in \cite{LibMar}, this defines a unique Poisson structure $\Lambda$ on $\dualJ$. In coordinates, we write
\[
\Lambda(dF,dG) = \{F(t,q,p),G(t,q,p)\} = \fpd{F}{q^i}\fpd{G}{p_i} - \fpd{F}{p_i}\fpd{G}{q^i},
\]
i.e.\ our sign convention is such that
\begin{equation}
\Lambda = \fpd{}{q^i}\wedge\fpd{}{p_i}. \label{Lambda}
\end{equation}
The corresponding Poisson map $P:\oneforms{\dualJ}\rightarrow\vectorfields{\dualJ}$ is defined by $\Lambda(\alpha,\beta)=\langle P(\alpha),\beta\rangle$, and another part of the sign convention we choose for is to put
\begin{equation}
X_F = - P(dF) = \fpd{F}{p_i}\fpd{}{q^i} - \fpd{F}{q^i}\fpd{}{p_i}. \label{XF}
\end{equation}
It is of interest to characterize $P$ also by its action on the local basis of 1-forms on $\dualJ$ which we used before. For $\alpha\in\oneforms{E}$ and $X\in\Vvectors$, we have
\begin{align}
P(\pi^*\alpha) &= \alpha_i\fpd{}{p_i} = \V{\alpha}, \label{P1} \\
P(dF_X) &= p_j\fpd{X^j}{q^i}\fpd{}{p_i} - X^i\fpd{}{q^i} = - \wt{X}, \label{P2}
\end{align}
and note for completeness that for the horizontal lift (\ref{hR}) of a $(1,1)$ tensor,
\begin{equation}
P(\h{R}) = p_iR^i_j\fpd{}{p_j} = \V{R}. \label{P3}
\end{equation}

Now let $R$ as before be a $(1,1)$ tensor on $E$ with the property $R(dt)=0$. We wish to investigate under what circumstances the complete lift $\wt{R}$ is a candidate to become a recursion operator for $P$, or expressed differently, for the couple $(P,\wt{R})$ to define a Poisson-Nijenhuis structure on $\dualJ$. For general aspects of Poisson-Nijenhuis structures see e.g.\ \cite{MaMor} and \cite{KosMa}. A preliminary condition to be satisfied is that $\wt{R}$ should commute with the Poisson map $P$. Now, using (\ref{tRa}), (\ref{P1}) and (\ref{tRv}), we get
\[
P\wt{R}(\pi^*\alpha) = P(\pi^*R(\alpha)) = \V{(R(\alpha))} = \wt{R}(\V{\alpha}) = \wt{R}P(\pi^*\alpha).
\]
Likewise, using (\ref{tRdF}), (\ref{P2}), (\ref{P3}) and (\ref{tRtX}), we have
\[
P\wt{R}(dF_X) = P\big(dF_{R(X)} - \h{(\lie{X}R)}\big) = - \wt{R(X)} - \V{(\lie{X}R)} = - \wt{R}(\wt{X}) = \wt{R}P(dF_X).
\]
This confirms the required commutation property, which is a condition also for the so-called Magri-Morosi concomitant $\mu_{\wt{R},P}$ \cite{MaMor} to be a tensor field of type $(1,2)$. Its defining relation then can be formulated as follows: for all $\sigma\in\oneforms{\dualJ}$ and $Z\in\vectorfields{\dualJ}$,
\begin{equation}
\mu_{\wt{R},P}(\sigma,Z) := \big(\lie{P(\sigma)}\wt{R}\big)(Z) - P\big(\lie{X}(\wt{R}(\sigma))\big) + P\big(\lie{\wt{R}(Z)}\sigma\big) \label{MM}
\end{equation}
is a vector field on $\dualJ$. The couple $(P,\wt{R})$ will define a Poisson-Nijenhuis structure if $\wt{R}$ has vanishing Nijenhuis torsion and $\mu_{\wt{R},P}=0$. We know all about vanishing $N_{\wt{R}}$. To check whether the Magri-Morosi concomitant vanishes, we further need the following list of properties of Lie derivatives of the 1-forms we obtained on $\dualJ$ by lifting operations.

\begin{lemma} Let $\alpha,\beta\in\oneforms{E}$, $X\in\Vvectors$, $Y\in\Vvectors\cup\tvectors$, while $R$ and $Q$ denote $(1,1)$ tensors on $E$ vanishing on $dt$. Then,
\begin{alignat}{3}
&\lie{\V{\beta}}(\pi^*\alpha) = 0 &\qquad &\lie{\V{Q}}(\pi^*\alpha) =0 &\qquad
&\lie{\wt{Y}}(\pi^*\alpha) = \pi^*(\lie{Y}\alpha) \label{nl1} \\
&\lie{\V{\beta}}dF_X = \pi^*di_X\beta &\quad &\lie{\V{Q}}dF_X = dF_{QX} &\quad
&\lie{\wt{Y}}dF_X = dF_{[Y,X]} \label{nl2} \\
&\lie{\V{\beta}} \h{R} = \pi^*R(\beta) &\quad &\lie{\V{Q}} \h{R} = \h{(Q\circ R)} &\quad
&\lie{\wt{Y}} \h{R} = \h{(\lie{Y}R)}. \label{nl3}
\end{alignat}
\end{lemma}
\proof\ The proof is a straightforward calculation which can be done either by evaluating both sides on the usual basis of vector fields and making use of Lie derivative properties obtained before, or perhaps more simply by a direct coordinate calculation. \qed

\begin{prop} The Magri-Morosi concomitant $\mu_{\wt{R},P}$ vanishes identically.
\end{prop}
\proof\ For an elegant proof, we simply let $\sigma$ and $Z$ in the defining relation (\ref{MM}) run over the set of 1-forms and vector fields on $\dualJ$ which we have used all the time to generate a local basis. For $\sigma=\pi^*\alpha$ and $Z=\V{\beta}$, with $\alpha,\beta\in\oneforms{E}$, one easily checks that all three terms vanish separately. For $\sigma=\pi^*\alpha$ and $Z=\wt{Y}$, using (\ref{aux}) and (\ref{nl1}b,\ref{nl1}c), the right-hand side of (\ref{MM}) reduces in the first place to
\[
-\V{(\lie{RY}\alpha)} + \V{(\lie{Y}(R\alpha))} - P(\pi^*\lie{Y}(R\alpha)) + P(\pi^*\lie{RY}\alpha),
\]
and this is clearly zero in view of (\ref{P1}). Next, for $\sigma=dF_X$ and $Z=\V{\beta}$, using (\ref{P2}) in the first term, (\ref{tRdF}) in the second and (\ref{tRv}) in the third, we get
\[
-\big(\lie{\wt{X}}\wt{R}\big)(\V{\beta}) - P\big(\lie{\V{\beta}}(dF_{RX} - \h{(\lie{X}R)})\big) + P\big(\lie{\V{(R\beta)}}dF_X\big).
\]
We subsequently use (\ref{lieR1}) in the first term, (\ref{nl2}a) and (\ref{nl3}a) in the second and (\ref{nl2}a) in the third again, and when we next evaluate the $P$-terms everything cancels out again. Finally, after a first evaluation, we get for $\mu_{\wt{R},P}(dF_X,\wt{Y})$:
\[
-\big(\lie{\wt{X}}\wt{R}\big)(\wt{Y}) - P\big(\lie{\wt{Y}}(dF_{RX} - \h{(\lie{X}R)})\big) + P\big(\lie{\wt{RY}}dF_X + \lie{\V{(\lie{Y}R)}}dF_X\big).
\]
It is then a matter of using (\ref{nl2}c), (\ref{nl3}c) and (\ref{nl2}b), plus the properties (\ref{P2}) and (\ref{P3}) of $P$, to come to an expression where everything cancels out again. \qed

\begin{thm} Let $R$ be a $(1,1)$ tensor on $\tau:E\rightarrow \R$ with the property $R(dt)=0$ and let $\wt{R}$ be its complete lift to $\dualJ$. Denote by $P:\oneforms{\dualJ}\rightarrow \vectorfields{\dualJ}$ the canonical Poisson map on $\dualJ$. Then, $(P,\wt{R})$ is a Poisson-Nijenhuis structure on $\dualJ$ if and only if $N_R=0$, where $N_R$ is the Nijenhuis torsion of $R$ on $E$.
\end{thm}
\proof\ We have shown that $P\wt{R}=\wt{R}P$ and that the Magri-Morosi concomitant $\mu_{\wt{R},P}$ vanishes. The only other requirement for having a Poisson-Nijenhuis structure then is that $N_{\wt{R}}=0$. But Theorem~2 tells us that this is equivalent to $N_R=0$. \qed

\section{Darboux-Nijenhuis coordinates}

It is known that under appropriate conditions, there exist special coordinates which are simultaneously adapted to the Poisson structure and the recursion operator of a Poisson-Nijenhuis structure, in the sense that they diagonalize the recursion tensor and provide canonical (Darboux) coordinates for the Poisson tensor. A few general references in this respect are \cite{Magri} and \cite{FP2}. Yet it is hard to find a detailed explanation on the way such coordinates will arise. So we will try to give such details here for the case of our $(P,\wt{R})$ structure.

The basic assumption is that the tensor $R$ on $E$ is algebraically diagonalizable with distinct eigenvalues. Observe then that the property $R(dt)=0$ immediately says that $dt$ is an eigenform corresponding to the eigenvalue $\lambda_0=0$. Obviously, $R$ is degenerate, but our assumption implies that it has rank $n$. We write $(R^\alpha_\beta)$ now for the matrix representation of $R$ (as linear map on vectors), with Greek indices running from 0 to $n$, and $\alpha$ in the role of row index. So, with $R$ as in (\ref{R}), the first row of the matrix has only zeros, the $R^i_0$ constitute the remaining elements of the first column, and the $n\times n$ matrix $(R^i_j)$ is non-singular. Since the Poisson tensor $\Lambda$ already takes its canonical form (\ref{Lambda}) in natural bundle coordinates $(t,q,p)$ on $\dualJ$, what we are after is a coordinate transformation which does not destroy this canonical form and achieves the diagonalization of $\wt{R}$ in coordinates. Hence, i
 t will h
 ave to be the induced transformation of a time-dependent coordinate transformation on $E$ which diagonalizes $R$. From the fundamental paper of Fr\"olicher and Nijenhuis \cite{FroNij}, we learn that diagonalizability in coordinates requires vanishing of the so-called Haantjes tensor. But a direct application of this theory in our case, where the manifold $E$ has coordinates $x^\alpha=(t,q^i)$ say, will merely guarantee the existence of new coordinates $y^\beta=y^\beta(x^\alpha)$ which do the job. This is a supplementary reason for going through the procedure in some detail here, because we need a transformation which preserves the fibred structure of $E$, i.e.\ the $y^\beta$ should be of the form $(t,Q^i(t,q))$.

Let $X_{(\alpha)}$ be a local basis for $\vectorfields{E}$ consisting of eigenvectors of $R$ (the extra brackets used for the index are meant to indicate that there are no summations over repeated indices in what follows). Then,
\begin{align*}
N_R(X_{(\alpha)},X_{(\beta)}) &= (R-\lambda_{(\alpha)})(R-\lambda_{(\beta)})([X_{(\alpha)},X_{(\beta)}])  \\
& \hspace*{3mm} + (\lambda_{(\alpha)}-\lambda_{(\beta)}) \big(X_{(\alpha)}(\lambda_{(\beta)})\,X_{(\beta)} + X_{(\beta)}(\lambda_{(\alpha)})\,X_{(\alpha)}\big).
\end{align*}
Following \cite{FroNij}, we next look at the Haantjes tensor, defined by
\[
\mathcal{H}_R(X,Y) := R^2N_R(X,Y) + N_R(RX,RY) - RN_R(RX,Y) - RN_R(X,RY),
\]
and easily obtain that
\begin{align*}
\mathcal{H}_R (X_{(\alpha)},X_{(\beta)}) &= (R-\lambda_{(\alpha)})(R-\lambda_{(\beta)})N_R(X_{(\alpha)},X_{(\beta)}) \\
&= (R-\lambda_{(\alpha)})^2(R-\lambda_{(\beta)})^2([X_{(\alpha)},X_{(\beta)}]).
\end{align*}
Since $N_R=0$, also $\mathcal{H}_R=0$ and this is the necessary and sufficient condition for a diagonalizable $R$ to be diagonalizable in coordinates. Indeed, we see that in our case $\mathcal{H}_R=0$ implies that $[X_{(\alpha)},X_{(\beta)}]
\in \sp \{X_{(\alpha)},X_{(\beta)}\}$ for all $\alpha,\beta$. If we denote by $\D{\alpha}$ the distribution spanned by $X_{(\alpha)}$, Frobenius theorem implies in such a case that all $\D{\alpha}$ are simultaneously integrable. In other words, there exist new coordinates $y^\alpha$ such that in those coordinates, $\D{\alpha}=\sp \{\partial/\partial {y^\alpha}\}$. Our extra concern now is that such new coordinates should be of the form indicated above. To see that this is possible, it suffices to look at the dual picture of eigenforms. Let $\D{\alpha}^\perp$ denote the annihilator of $\D{\alpha}$ and put $\D{\alpha}^*= \cap_{\beta\neq\alpha}\D{\beta}^\perp$. Then by construction $\langle X_{(\beta)}, \rho_{(\alpha)}\rangle =0\ \forall \beta\neq\alpha$ and $\rho_{(\alpha)}\in\D{\alpha}^*$, while
\[
\langle R(X_{(\alpha)}),\rho_{(\alpha)}\rangle = \lambda_{(\alpha)}\langle X_{(\alpha)},\rho_{(\alpha)}\rangle = \langle X_{(\alpha)},R(\rho_{(\alpha)})\rangle.
\]
Hence, $\rho_{(\alpha)}$ is an eigenform of $R$ corresponding to the eigenvalue $\lambda_{(\alpha)}$. Therefore, in the coordinates $y^\alpha$ simultaneously adapted to all $\D{\alpha}$, $\rho_{(\alpha)}$ is in the module generated by $dy^\alpha$. So, in the dual picture, the coordinate transformation $(t,q^i)\rightarrow y^\alpha$ has the task of producing eigenforms of the form $dy^\alpha$. But we know that $dt$ is an eigenform (with eigenvalue zero), hence we can simply take $y^0=t$, meaning that $R$ will indeed be diagonalized by a transformation of the form $(t,q)\rightarrow (t,Q(t,q))$. In the new variables, the eigenvectors can be taken to be coordinate fields, so that they commute and the Nijenhuis tensor expression reduces to
\[
N_R(X_{(\alpha)},X_{(\beta)}) = (\lambda_{(\alpha)}-\lambda_{(\beta)}) \big(X_{(\alpha)}(\lambda_{(\beta)})\,X_{(\beta)} + X_{(\beta)}(\lambda_{(\alpha)})\,X_{(\alpha)}\big).
\]
It follows that the stronger property $N_R=0$ now implies that $X_{(\alpha)}(\lambda_{(\beta)})=0$ for $\alpha\neq\beta$. The conclusion is that each eigenvalue $\lambda_{(\alpha)}$ is a function of the corresponding coordinate $y^\alpha$ only. The final conclusion is that in the new coordinates $(t,Q)$, $R$ will take the form
\[
R= \sum_{i=1}^n \lambda_{(i)}(Q^i) \fpd{}{Q^i} \otimes dQ^i,
\]
and this implies that
\[
\wt{R} = \sum_{i=1}^n \lambda_{(i)}(Q^i) \left(\fpd{}{Q^i} \otimes dQ^i + \fpd{}{P_i}\otimes dP_i \right) .
\]
As said before, the nature of the coordinate transformation involved in this process ensures that the Poisson tensor $\Lambda$ will still have the canonical form
\[
\Lambda = \fpd{}{Q^i}\wedge \fpd{}{P_i},
\]
so we have indeed obtained Darboux-Nijenhuis coordinates for the Poisson-Nijenhuis structure $(P,\wt{R})$.

\section{Concluding remarks}

Poisson-Nijenhuis structures play a prominent role in the study and characterization of integrable Hamiltonian systems, bi-Hamiltonian systems and Hamilton-Jacobi separability. Some interesting general references in this respect are \cite{NunezMarle} and \cite{FP}. More specifically, when it comes to aspects of standard Hamilton-Jacobi separability for autonomous Hamiltonian systems, the recursion operator of a Poisson-Nijenhuis structure will generally be the complete lift to the cotangent bundle $T^*M$ of a type $(1,1)$ tensor on $M$. An interesting subcase of St\"ackel separability arises when the $(1,1)$ tensor on $M$ is a so-called special conformal Killing tensor with respect to the Riemannian metric determining the kinetic energy of the system. Such tensors, also called or at least closely related to what are called Benenti tensors (after Benenti's pioneering work in \cite{Ben} and \cite{Ben2}), make their appearance, for example, in \cite{IMM} and \cite{Crampin}. They
  automat
 ically satisfy the condition of vanishing Nijenhuis torsion. Special conformal Killing tensors further play a significant role in related work, such as the study of a certain bi-differential calculus \cite{CST}, and the intrinsic characterization and generalization in \cite{CraSar1} of what Lundmark called Newtonian systems of cofactor type (see \cite{Lundmark} and \cite{Lundmark2}). We recently also succeeded in providing a full geometrical description of so-called driven cofactor systems in \cite{SW}.

It is our intention to study in forthcoming papers various of the above mentioned aspects in the context of time-dependent Hamiltonian systems. For these objectives, understanding the notion of complete lift of a $(1,1)$ tensor from $E$ to $\dualJ$, and having extensive knowledge of its properties, is a necessary prerequisite which we hope to have achieved now.

\subsubsection*{Acknowledgements} This work is part of the IRSES
project GEOMECH (nr. 246981) within the 7th European Community
Framework Programme.

\end{document}